\pdfoutput=1
\documentclass{amsart}
\usepackage{float}
\usepackage[section]{placeins}
\restylefloat{table}

\usepackage{amsmath,amssymb,amsthm}
\usepackage{graphicx,tikz}
\usepackage{subcaption}

\newtheorem {theorem}{Theorem}
\newtheorem*{thm}{Theorem}
\newtheorem*{proposition}{Proposition}

\theoremstyle{definition}

\theoremstyle{remark}

\begin{document}

\title[]{Sequences of Well-distributed vertices on Graphs and Spectral Bounds on Optimal Transport }
\thanks{This paper is part of the author's PhD thesis at Yale University and has been partially supported by NSF (DMS-1763179).}

\author[]{Louis Brown}
\address{Department of Mathematics, Yale University, New Haven, CT 06511, USA}
\email{louis.brown@yale.edu}

\begin{abstract} Given a graph $G=(V,E)$, suppose we are interested in selecting a sequence of
vertices $(x_j)_{j=1}^n$ such that $\left\{x_1, \dots, x_k\right\}$ is `well-distributed' \textit{uniformly in $k$}.
We describe a greedy algorithm motivated by potential theory and corresponding developments in the continuous setting.
The algorithm performs nicely on graphs and may be of use for sampling problems. We can interpret the algorithm as trying to greedily minimize a negative Sobolev norm; we explain why this is related to Wasserstein distance by establishing a purely spectral bound on the Wasserstein
distance on graphs that mirrors R. Peyre's estimate in the continuous setting. We illustrate this with many examples and discuss several open problems.
\end{abstract}

\maketitle

\vspace{-10pt}

\section{Introduction}

\subsection{Introduction.}  The purpose of this paper is to discuss a basic problem on finite graphs $G=(V,E)$ that is already difficult on the unit interval $[0,1]$. Given a
metric space (say, the unit interval, a compact manifold or a finite graph), how does one construct a sequence $x_1, x_2, \dots$ of elements in the space such that their distribution is \textit{uniformly good}---by this
we mean that if one takes the first $k$ elements $\left\{x_1, \dots, x_k\right\}$, then this set is very nearly as evenly distributed on the set as any
set of $k$ elements would be. We have not clearly defined what notion of `well-distributed' we mean, this will depend on the actual setting; the question is frequently interesting for several different such notions.
There is not even a canonical answer on the unit interval $[0,1]$ but the van der Corput sequence 
$$ \frac{1}{2}, \frac{1}{4}, \frac{3}{4}, \frac{1}{8}, \frac{5}{8}, \frac{3}{8}, \frac{7}{8}, \dots$$
constructed via an inverse binary digit expansion (see \cite{dick}), is a good example.

\begin{center}
\begin{figure}[h!]
\begin{tikzpicture}[scale=1.2]
\draw [very thick] (0,0) -- (10,0);
\draw [very thick] (0,-0.1) -- (0,0.1);
\draw [very thick] (10,-0.1) -- (10,0.1);
\filldraw (5,0) circle (0.04cm);
\node at (5, 0.4) {$\frac{1}{2}$};
\node at (5, -0.4) {1};
\filldraw (2.5,0) circle (0.04cm);
\node at (2.5, 0.4) {$\frac{1}{4}$};
\node at (2.5, -0.4) {2};
\filldraw (7.5,0) circle (0.04cm);
\node at (7.5, 0.4) {$\frac{3}{4}$};
\node at (7.5, -0.4) {3};
\filldraw (10/8,0) circle (0.04cm);
\node at (10/8, 0.4) {$\frac{1}{8}$};
\node at (10/8, -0.4) {4};
\filldraw (5*10/8,0) circle (0.04cm);
\node at (5*10/8, 0.4) {$\frac{5}{8}$};
\node at (5*10/8, -0.4) {5};
\filldraw (3*10/8,0) circle (0.04cm);
\node at (3*10/8, 0.4) {$\frac{3}{8}$};
\node at (3*10/8, -0.4) {6};
\filldraw (7*10/8,0) circle (0.04cm);
\node at (7*10/8, 0.4) {$\frac{7}{8}$};
\node at (7*10/8, -0.4) {7};
\end{tikzpicture}
\caption{The first 7 elements of the van der Corput sequence.}
\end{figure}
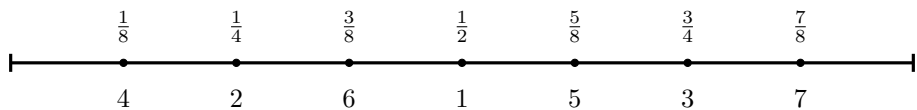
\end{center}
\vspace{-10pt}
We see that the first $k$ elements of the sequence are not as uniformly distributed as $k$ equispaced points but they are always fairly uniformly distributed independently of what $k$ is.  There are many different reasons why one could be interested in such sequences: they are natural sampling
points for functions (especially for on-line selection and in cases where one does not know in advance in how
many points one can sample) but there is also an obvious combinatorial question (`How well distributed can
sequences be? What is the unavoidable degree of irregularity?').

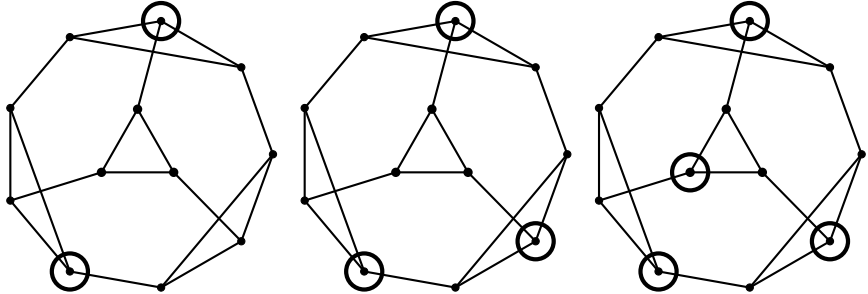
\begin{figure}[h!]
\begin{minipage}[l]{.3\textwidth}
\centering
  \begin{tikzpicture}[scale=0.6]
\foreach \a in {1,2,...,9}{
\filldraw (\a*360/9: 3cm) circle (0.08cm);
};
\foreach \a in {1,2,...,9}{
\draw [thick] (\a*360/9: 3cm) --  (\a*360/9 + 360/9: 3cm);
};
\filldraw [thick] (0.8, -0.4) circle (0.08cm);
\filldraw [thick] (-0.8, -0.4) circle (0.08cm);
\filldraw [thick] (0, 1) circle (0.08cm);
\draw [thick] (0.8, -0.4) --  (8*360/9: 3cm);
\draw [thick] (0.8, -0.4) -- (-0.8, -0.4) -- (0,1) -- (0.8, -0.4);
\draw [thick] (-0.8, -0.4) --  (5*360/9: 3cm);
\draw [thick] (0,1) --  (2*360/9: 3cm);
\draw [thick] (1*360/9: 3cm) --  (3*360/9: 3cm);
\draw [thick] (4*360/9: 3cm) --  (6*360/9: 3cm);
\draw [thick] (7*360/9: 3cm) --  (9*360/9: 3cm);
\draw [ultra thick]  (2*360/9: 3cm) circle (0.4cm);
\draw [ultra thick]  (6*360/9: 3cm) circle (0.4cm);
   \end{tikzpicture}
\end{minipage} 
\begin{minipage}[l]{.3\textwidth}
\centering
  \begin{tikzpicture}[scale=0.6]
\foreach \a in {1,2,...,9}{
\filldraw (\a*360/9: 3cm) circle (0.08cm);
};
\foreach \a in {1,2,...,9}{
\draw [thick] (\a*360/9: 3cm) --  (\a*360/9 + 360/9: 3cm);
};
\filldraw [thick] (0.8, -0.4) circle (0.08cm);
\filldraw [thick] (-0.8, -0.4) circle (0.08cm);
\filldraw [thick] (0, 1) circle (0.08cm);
\draw [thick] (0.8, -0.4) --  (8*360/9: 3cm);
\draw [thick] (0.8, -0.4) -- (-0.8, -0.4) -- (0,1) -- (0.8, -0.4);
\draw [thick] (-0.8, -0.4) --  (5*360/9: 3cm);
\draw [thick] (0,1) --  (2*360/9: 3cm);
\draw [thick] (1*360/9: 3cm) --  (3*360/9: 3cm);
\draw [thick] (4*360/9: 3cm) --  (6*360/9: 3cm);
\draw [thick] (7*360/9: 3cm) --  (9*360/9: 3cm);
\draw [ultra thick]  (2*360/9: 3cm) circle (0.4cm);
\draw [ultra thick]  (6*360/9: 3cm) circle (0.4cm);
\draw [ultra thick]  (8*360/9: 3cm) circle (0.4cm);
   \end{tikzpicture}
\end{minipage} 
\begin{minipage}[l]{.3\textwidth}
\centering
  \begin{tikzpicture}[scale=0.6]
\foreach \a in {1,2,...,9}{
\filldraw (\a*360/9: 3cm) circle (0.08cm);
};
\foreach \a in {1,2,...,9}{
\draw [thick] (\a*360/9: 3cm) --  (\a*360/9 + 360/9: 3cm);
};
\filldraw [thick] (0.8, -0.4) circle (0.08cm);
\filldraw [thick] (-0.8, -0.4) circle (0.08cm);
\filldraw [thick] (0, 1) circle (0.08cm);
\draw [thick] (0.8, -0.4) --  (8*360/9: 3cm);
\draw [thick] (0.8, -0.4) -- (-0.8, -0.4) -- (0,1) -- (0.8, -0.4);
\draw [thick] (-0.8, -0.4) --  (5*360/9: 3cm);
\draw [thick] (0,1) --  (2*360/9: 3cm);
\draw [thick] (1*360/9: 3cm) --  (3*360/9: 3cm);
\draw [thick] (4*360/9: 3cm) --  (6*360/9: 3cm);
\draw [thick] (7*360/9: 3cm) --  (9*360/9: 3cm);
\draw [ultra thick]  (2*360/9: 3cm) circle (0.4cm);
\draw [ultra thick]  (6*360/9: 3cm) circle (0.4cm);
\draw [ultra thick]  (8*360/9: 3cm) circle (0.4cm);
\draw [ultra thick]  (-0.8, -0.4) circle (0.4cm);
   \end{tikzpicture}
\end{minipage} 
\caption{The first $k$ (here $k=2,3,4$) elements of the sequence are nearly as evenly distributed as any set of $k$ vertices could be.}
\end{figure}

Historically, the question has lead to remarkable connections to other fields of mathematics such as Number Theory, Harmonic Analysis and even Probability Theory (we will discuss several of these connections in \S 2).  We will now state one informal version of the main problem before stating a precise version further below.

\begin{quote} \textbf{Main Problem} (informal version)\textbf{.}  Given a finite graph $G=(V,E)$, how
would one select a sequence of vertices that are \textit{uniformly good}?  In what metric would one
measure the `goodness' of such a sequence?
\end{quote}
\begin{center}
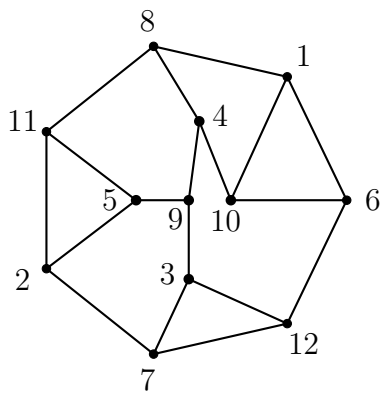
\begin{figure}[h!]
  \begin{tikzpicture}[scale=0.7]
\foreach \a in {1,2,...,7}{
\filldraw (\a*360/7: 3cm) circle (0.08cm);
};
\foreach \a in {1,2,...,7}{
\draw [thick] (\a*360/7: 3cm) --  (\a*360/7 + 360/7: 3cm);
};
\filldraw [thick] (0.8, 0) circle (0.08cm);
\draw [thick] (0.8, 0) --  (1*360/7: 3cm);
\draw [thick] (0.8, 0) --  (7*360/7: 3cm);
\draw [thick] (0.8, 0) --  (0.2, 1.5);
\filldraw [thick] (0.2, 1.5) circle (0.08cm);
\filldraw [thick] (0,0) circle (0.08cm);
\filldraw [ thick] (-1,0) circle (0.08cm);
\filldraw [ thick] (0,-1.5) circle (0.08cm);
\draw [ thick] (0.2, 1.5) --  (2*360/7: 3cm);
\draw [ thick] (0.2, 1.5) --  (0,0);
\draw [ thick] (-1, 0) --  (0,0);
\draw [ thick] (0, -1.5) --  (0,0);
\draw [ thick] (-1, 0) --  (3*360/7: 3cm);
\draw [ thick] (-1, 0) --  (4*360/7: 3cm);
\draw [ thick] (0,-1.5) --  (5*360/7: 3cm);
\draw [ thick] (0,-1.5) --  (6*360/7: 3cm);
\node at (1*360/7: 3.5cm) {\Large 1};
\node at (4*360/7: 3.5cm) {\Large 2};
\node at (-0.4,-1.4) {\Large 3};
\node at (0.6,1.6) {\Large 4};
\node at (-1.5,0) {\Large 5};
\node at (7*360/7: 3.5cm) {\Large 6};
\node at (5*360/7: 3.5cm) {\Large 7};
\node at (2*360/7: 3.5cm) {\Large 8};
\node at (-0.25,-0.35) {\Large 9};
\node at (0.7,-0.4) {\Large 10};
\node at (3*360/7: 3.5cm) {\Large 11};
\node at (6*360/7: 3.5cm) {\Large 12};
   \end{tikzpicture}
\caption{The Frucht Graph (see \S4.3) and the enumeration obtained by the algorithm when starting with vertex 1.}
\end{figure}
\end{center}
Fig. 2 contains a simple such example: taking the Truncated Tetrahedral Graph (see \S4.2), in which order should one select the vertices so as to obtain a sequence that is uniformly evenly distributed? Fig. 3 showcases another example, on the Frucht Graph (see \S4.3).  In Fig. 5, we display an example on the Nauru Graph, a $3$-regular bipartite graph with 24 vertices and 36 edges.  Even before making the notion of quality precise, we can get some intuition from these simple examples. The enumeration of the vertices in each case was generated by the algorithm discussed below, in \S2.2.  The organization of this article is as follows: in \S 1.2 we introduce the Wasserstein distance, and state Kantorovich-Rubinstein duality in the setting of finite graphs.  In \S2, we introduce a novel algorithm for selecting evenly distributed vertices on graphs, and in \S2.3 we provide a Theorem quantitatively demonstrating the even distribution of these vertices.  In \S3 we explore the connection between optimal transport cost and spectral properties of the Laplacian, providing in \S3.2 a Theorem bounding the former in terms of the latter on graphs.  In \S3.4 and 3.5 we conduct a case study of cycle graphs and torus grid graphs, respectively, in the context of this Theorem.  In \S4, we provide numerics displaying the performance of the algorithm on a variety of different graphs.  In \S5, we present proofs of both Theorems.  In \S6, we list connections to results in the continuous setting.\vspace{-1 em}
\subsection{Wasserstein Distance.}
\begin{figure}[H]\begin{tikzpicture}[scale=.7]
\foreach \a in {1,...,6}{
\filldraw (\a*60:3cm) circle (.08cm);
\node at (\a*60:3.6cm){\Large{$x_\a$}};
\draw [thick] (\a*60+60:3cm)--(\a*60:3cm);}
\draw [ultra thick] (60:3cm) circle (.3cm);
\draw [ultra thick] (240:3cm) circle (.3cm);
\node at (60:2.3cm) {$1/2$};
\node at (240:2.3cm) {$1/2$};
\end{tikzpicture}
\caption{$W_1((\delta_{x_1}+\delta_{x_2})/2,\sum_{k=1}^6\delta_{x_k}/6)=2/3$: We can transport $1/6$ units of mass from $x_1$ to each of $x_2$ and $x_6$, and similarly with $x_4$ to $x_3$ and $x_5$, incurring a total cost of $4\times1/6$.}\end{figure}
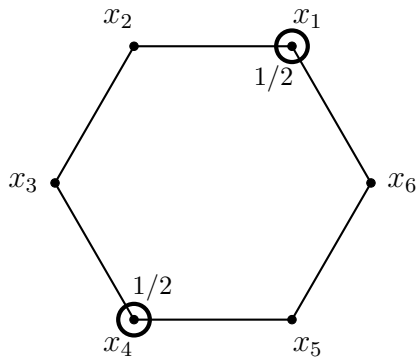

Wasserstein Distance is a notion of distance on probability distributions introduced by Wasserstein in 1969 \cite{wasser}. Its simplest instance is $W_1(\mu, \nu)$, also known as Earth Mover's Distance, which is defined as the minimal amount of cost required to transport one distribution $\mu$ to match another distribution $\nu$---here, cost is defined as mass $\times$ distance: transporting $\varepsilon$ units of mass over a distance 
of $\delta$ has a cost of $\varepsilon\delta$.  The notion of Wasserstein distance can be immediately applied to finite graphs: we will work only in the simple case of unweighted graphs where the distance between two vertices is given as the length of the shortest path connecting them; extensions to the weighted case are conceivable.  In short, transporting $\varepsilon$ mass over a single edge has a $W_1$ cost of $\varepsilon$ (see Fig. 4).
More formally, on an abstract metric space $X$ equipped with a metric $d$, we define
$$ W_p(\mu, \nu) = \left( \inf_{\gamma \in \Gamma(\mu, \nu)} \int_{X \times X}{d(x,y)^p d \gamma(x,y)}\right)^{1/p},$$
where $\Gamma(\mu, \nu)$ denotes the collection of all measures on $X \times X$
with marginals $\mu$ and $\nu$, respectively (also called the set of all couplings of $\mu$ and $\nu$). Note that $\mu,\nu$ need not be probability measures.  We will, throughout this paper, work exclusively with the Earth Mover's Distance $W_1$ (although extensions to more general $W_p$ are certainly conceivable). The Earth Mover's Distance is particularly nice to work with: by Kantorovich-Rubinstein duality (see e.g., \cite{villani})
$$ W_1(\mu, \nu) = \sup\left\{ \int_{X}{ f d\mu} - \int_{X}{ f d\nu}: f~\mbox{is 1-Lipschitz} \right\}.$$  Because of this, $W_1$ is translation invariant---for positive measures $\mu,\nu, \mu'$, we have $$W_1(\mu,\nu)=W_1(\mu+\mu',\nu+\mu').$$  Thus, if we have a positive measure $\mu$ decomposed into non-positive measures $\mu_1,\mu_2$ as $\mu=\mu_1+\mu_2$, we may use this translation invariance to write $$W_1(\mu,\nu)=W_1(\mu_1^++\mu_2^+,\nu+\mu_1^-+\mu_2^-),$$ where $\mu_i^+=\max\{\mu_i,0\}$ is the positive part of $\mu_i$ and $\mu_i^-=\max\{-\mu_i,0\}$ is the negative part.
Now that we have defined the Wasserstein distance, we may make precise our original question regarding evenly distributed vertices:
\begin{quote}  \textbf{Main Problem} (formal version)\textbf{.}   Given a finite graph $G=(V,E)$, how
would one select a sequence of vertices such that
$$ W_1\left( \frac{1}{k}\sum_{j=1}^k{\delta_{x_j}}, dx\right) \qquad \mbox{is small for all}~k,$$
where $dx$ is normalized counting measure having weight
$|V|^{-1}$ on each vertex of $G$.
\end{quote}
How small one could expect this quantity to be will depend on the particular geometry of the graph. For simplicity of exposition, consider $M=\mathbb T^d$, the $d$-dimensional torus with normalized volume measure $dx$: it is a basic exercise to show that the Earth Mover's Distance satisfies
$$ W_1\left( \frac{1}{n}\sum_{j=1}^{n}{\delta_{x_j}}, dx\right) \geq c_d n^{-1/d}$$
for all sets of points $\left\{x_1, \dots, x_n\right\} \subset M$. This clearly shows that the geometry (here: the dimension $d$) plays a role in what we can expect. We also emphasize that it is almost surely the case that our main question (as asked in \S 1.1) is of interest also for many other ways of making the notion of distribution quantitative and Wasserstein distance may be one of many (though certainly a rather canonical one).  We conclude our short introduction to Wasserstein distance by stating Kantorovich-Rubinstein duality, mentioned above, when $X$ is a finite graph.
\begin{proposition}[Kantorovich-Rubinstein, see e.g., \cite{kant,peyre}] Let $G=(V,E)$ be a finite, simple graph, let $f:V \rightarrow \mathbb{R}$ and let $W \subset V$ be a subset of vertices. Then
$$ \left| \frac{1}{|V|}\sum_{x\in V}{f(x)} - \frac{1}{|W|}\sum_{x\in W}{f(x)} \right| \leq W_1\left( \frac{1}{|W|}\sum_{x\in W}{\delta_{x}}, dx\right) \max_{x_i\sim x_j}{  |f(x_i) - f(x_j)|}.$$
\end{proposition}
This shows that our notion of uniform distribution of a subset of vertices has a natural connection to the question of sampling on graphs (i.e., reconstructing the average value of a `smooth' function by sampling in a subset of the vertices). The theory of sampling on graphs is in its infancy but rapidly developing, we refer to \cite{hu, linderman, pes1, pes2, shu}.

\begin{center}
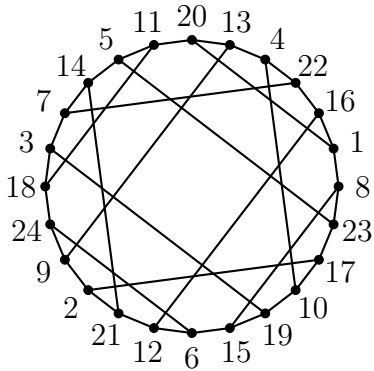
\begin{figure}[H]
  \begin{tikzpicture}[scale=0.65]
\foreach \a in {1,2,...,24}{
\filldraw (\a*360/24: 3cm) circle (0.09cm);
};
\foreach \a in {1,2,...,24}{
\draw [thick] (\a*360/24: 3cm) --  (\a*360/24 + 360/24: 3cm);
};
\draw [thick] (1*360/24: 3cm) -- (6*360/24: 3cm);
\draw [thick] (2*360/24: 3cm) -- (17*360/24: 3cm);
\draw [thick] (3*360/24: 3cm) -- (10*360/24: 3cm);
\draw [thick] (4*360/24: 3cm) -- (21*360/24: 3cm);
\draw [thick] (5*360/24: 3cm) -- (14*360/24: 3cm);
\draw [thick] (7*360/24: 3cm) -- (12*360/24: 3cm);
\draw [thick] (8*360/24: 3cm) -- (23*360/24: 3cm);
\draw [thick] (9*360/24: 3cm) -- (16*360/24: 3cm);
\draw [thick] (11*360/24: 3cm) -- (20*360/24: 3cm);
\draw [thick] (13*360/24: 3cm) -- (18*360/24: 3cm);
\draw [thick] (15*360/24: 3cm) -- (22*360/24: 3cm);
\draw [thick] (19*360/24: 3cm) -- (24*360/24: 3cm);
\node at (1*360/24: 3.5cm) {\Large 1};
\node at (15*360/24: 3.5cm) {\Large 2};
\node at (11*360/24: 3.5cm) {\Large 3};
\node at (4*360/24: 3.5cm) {\Large 4};
\node at (8*360/24: 3.5cm) {\Large 5};
\node at (18*360/24: 3.5cm) {\Large 6};
\node at (10*360/24: 3.5cm) {\Large 7};
\node at (24*360/24: 3.5cm) {\Large 8};
\node at (14*360/24: 3.5cm) {\Large 9};
\node at (21*360/24: 3.5cm) {\Large 10};
\node at (7*360/24: 3.5cm) {\Large 11};
\node at (17*360/24: 3.5cm) {\Large 12};
\node at (5*360/24: 3.5cm) {\Large 13};
\node at (9*360/24: 3.5cm) {\Large 14};
\node at (19*360/24: 3.5cm) {\Large 15};
\node at (2*360/24: 3.5cm) {\Large 16};
\node at (22*360/24: 3.5cm) {\Large 17};
\node at (12*360/24: 3.5cm) {\Large 18};
\node at (20*360/24: 3.5cm) {\Large 19};
\node at (6*360/24: 3.5cm) {\Large 20};
\node at (16*360/24: 3.5cm) {\Large 21};
\node at (3*360/24: 3.5cm) {\Large 22};
\node at (23*360/24: 3.5cm) {\Large 23};
\node at (13*360/24: 3.5cm) {\Large 24};
   \end{tikzpicture}
\caption{The Nauru Graph on 24 vertices: algorithm starting at 1.}
\end{figure}
\end{center}
\vspace{-20pt}
\section{The Algorithm}

\subsection{Setup}

We recall that for a finite graph $G=(V,E)$, we can define the adjacency matrix 
$$ A = \left(a_{ij}\right)_{i,j=1}^{|V|} \qquad \mbox{where}\qquad a_{ij}=\begin{cases} 1 \qquad &\mbox{if}~x_i \sim_{E} x_j \\ 0 \qquad &\mbox{otherwise} \end{cases}$$
as well as the degree matrix
$$ D = \left(d_{ij}\right)_{i,j=1}^{|V|} \qquad \mbox{where}\qquad d_{ij}=\begin{cases} \mbox{deg}(x_i) \qquad &\mbox{if}~i = j \\ 0 \qquad &\mbox{otherwise.} \end{cases}$$
With these definitions, we can define a notion of a Laplacian via
$$ L = D-A.$$
We denote the eigenvectors of $L$ by $\phi_i$ with corresponding eigenvalues $\lambda_i$, i.e., $L\phi_i=\lambda_i\phi_i$.  Note that $L$ has all its eigenvalues in $[0,2\max_{v\in V}\deg(v)]$. Since $L$'s columns sum to 0, we have for all measures $\mu$ that $L\mu$ has no net mass, i.e., is orthogonal to the constant vector, or has mean 0.  Since $L$ is symmetric, it is diagonalizable and all pairs of eigenvectors with distinct eigenvalues are orthogonal.  The induced ordering of eigenvectors is analogous to the continuous case:  small eigenvalue means slow oscillation frequency and the oscillation increases with the eigenvalue---the larger the eigenvalue, the more oscillation there is.  In particular, $\phi_1$ is constant with $\lambda_1=0$.  Since we assume $G$ is connected, there is only one instance of the trivial eigenvalue.  As $L$ is diagonalizable, we can also take arbitrary powers and define the fractional Laplacian $L^\alpha$ for $\alpha>0$.  Setting $n=|V|$, $$L^\alpha v=\sum_{i=1}^n \left\langle v,\phi_i\right\rangle\lambda_i^\alpha\phi_i.$$  To define $L^{-\alpha}$, we need to adjust this definition slightly: since $\lambda_1=0$, we first shift $v$ down by its mean to avoid dividing by 0.  (In other words, we simply ignore the component of $v$ in the direction of the constant eigenvector $\phi_1$.)  That is, $$L^{-\alpha} v=\sum_{i=2}^n\left\langle v,\phi_i\right\rangle\lambda_i^{-\alpha}\phi_i.$$  Of course, multiplying a vector by a matrix can be interpreted as applying an operator to a function, since vectors indexed by vertices are simply functions $V\to\mathbb R$.  Note that $AD^{-1}=I-LD^{-1}$ is a diffusion operator: each vertex splits its mass uniformly among its neighbors, and hence mass is preserved.  This operator has eigenvalues in $[-1,1]$.  If we instead apply the transpose $D^{-1}A$, this corresponds to each vertex \textit{taking} an equal portion of each of its neighbors masses, which, in general, is not mass-preserving. If $G$ is $k$-regular (each vertex has equal degree $k$), then $D=kI$ is scalar and these two notions coincide and equal $\frac{1}{k}A$.  We refer to \cite{chung, grig} for a good introduction to these notions and many references.
\subsection{Description of the Algorithm}
We present an algorithm, parametrized by $0<\alpha<1$, for greedily picking well-distributed vertices $x_k$ on graphs.  (If $\alpha\gg 1$ the algorithm degenerates and repeatedly selects the same small subset of distant vertices over and over.)  First, $x_1$ is chosen arbitrarily.  Then, vertices are picked recursively according to the following: $$ x_{k+1} = \arg\min_{x\in V} \left( L^{-\alpha} \sum_{j=1}^{k}{\delta_{x_j}}\right)(x),$$ breaking ties arbitrarily (Figures 2, 3, 5, and 9-11 display applications of this algorithm to various graphs).  If we write out the algorithm explicitly in terms of the spectrum of the Laplacian operator $L=D-A$, it becomes $$x_{k+1}=\arg\min_{x\in V}\sum_{i=2}^n \sum_{j=1}^k\frac{\phi_i(x_j)}{\lambda_i^{\alpha}}\phi_i(x),$$ where the $\phi_i$ are normalized with $\|\phi_i\|_{L^2}=1$, and we skip the constant eigenvector $\phi_1$ since it has eigenvalue 0 (note that this choice, while seemingly arbitrary, has no impact on the algorithm as any contribution from the constant vector can be ignored when computing $\arg\min$---the algorithm is independent of choice of right inverse of the Laplacian).  That is, we add up the projections of the indicator vector of the current vertex set, $\sum_{j=1}^k\delta_{x_j}$, onto each eigenvector, scaling down by the $\alpha$ power of the respective eigenvalue.  Consider the case of a cycle graph: if the number of vertices is sufficiently large, this is well approximated by a torus.  Setting $\alpha=1/2$ and identifying the torus with $[0,2\pi)/\sim:=\mathbb T$, we have a simple explicit formula for the inverse Laplacian of a point mass (see \S3.4 for the derivation): $$L^{-1/2}(\delta_{x_k})=-\frac{1}{\pi}\ln|2\sin((x-x_k)/2)|.$$  Note that $2\sin((x-x_k)/2)$ is precisely the Euclidean distance between the points at angles $x_i$ and $x$ on the unit circle (i.e., $|e^{2\pi ix}-e^{2\pi ix_k}|$).  Thus, the algorithm is simply maximizing the product of distances between points on the circle, by setting $$x_{k+1}=\arg\min_{x\in\mathbb T}\left(-\frac{1}{\pi}\sum_{j=1}^k \ln|2\sin((x-x_j)/2)|\right).$$  The arising sequence appears to behave on par with provably optimally regular sequences, and Steinerberger recently proved strong results on the regularity of such a sequence in \cite{steinii}, using techniques which are specific to this setting and unlikely to generalize to other graphs.  (We explore this example in more detail in $\S 3.4$.)  Nonetheless, these remarkable results on the torus and cycle graph give us hope that the algorithm may work comparably well on graphs more generally.

\subsection{A Theoretical Guarantee}

We prove a theoretical guarantee, for any finite graph $G$, that these sequences do exhibit at least a certain degree of regularity.

\begin{theorem}
Let $G$ be a simple connected graph and let $$\mu_k=\frac{1}{k}\sum_{j=1}^{k}\delta_{x_j},$$  where vertices $x_j$ are selected as above.  Then, for all $1\le k\le n$, $$\sum_{i=2}^n\frac{|\langle\mu_k,\phi_i\rangle|^2}{\lambda_i^{2\alpha}}\le \left(\max_{j\le k}\left\lVert L^{-2\alpha}\left(\delta_{x_j}\right)\right\rVert_{\ell^2}^2\right)k^{-1}.$$
\end{theorem}

\textit{Remark 1.}
Observe for the sake of comparison that  $$\sum_{i=2}^n \left|\langle\mu_k,\phi_i\rangle\right|^2=\|\mu_k\|^2_{\ell^2}-|\langle\mu_k,\phi_1\rangle|^2=\frac{1}{k}-\frac{1}{n}.$$  Thus, it is natural that the bound in the Theorem should be $\sim k^{-1}$.  However, $\lambda_i$ (and thus $\lambda_i^{2\alpha}$) may be arbitrarily close to 0, scaling up the terms in the sum substantially.  The only way to prevent this is for $\mu_k$ to be almost orthogonal to low-frequency eigenfunctions (which is cf. Erd\H{o}s-Tur\'an  \cite{erd1, erd2} a natural way of defining regularity, as it means that $\mu_k$ is concentrated at high frequencies).
\\\textit{Remark 2.}
Note that $$\max_{j\le k}\left\lVert L^{-2\alpha}\left(\delta_{x_j}\right)\right\rVert_{\ell^2}^2\le \max_{x\in V}\left\lVert L^{-2\alpha}\left(\delta_{x}\right)\right\rVert_{\ell^2}^2.$$ The term on the right side is an interesting quantity in itself, and there may be good bounds for it in terms of the geometry of the graph. As can be seen from the expansion into eigenfunctions, this quantity measures, implicitly, how much low-frequency eigenfunctions concentrate in a particular vertex.  In vertex-transitive graphs like cycle graphs and torus grid graphs we see that the quantity is actually independent of the vertex $x$.

\section{Spectral Bounds on Transport Distances}

\subsection{Motivation} 
The motivation behind the algorithm is two-fold:
\begin{enumerate}
\item The greedy algorithm tries to minimize a Sobolev norm
$$ \left\| \frac{1}{k} \sum_{j=1}^{k} \delta_{x_j} \right\|_{\dot H^{-1}}.$$ 
\item A recent result of R. Peyre \cite{pey} shows that, in the continuous setting,
$$ W_2(\mu,dx) \lesssim \|\mu\|_{\dot H^{-1}}.$$
\end{enumerate}
The purpose of this section is to establish a connection between problems of optimal transport and spectral properties of the Laplacian.
This is known to hold in the continuous case, we recall the following bound:
\begin{thm}[Carroll, Massaneda, Ortega-Cerda \cite{carr}] 
Let $(M,g)$ be a compact Riemannian manifold with normalized volume measure $dx$ and $\partial M = \emptyset$. 
 If $-\Delta_g \phi=\lambda \phi$ on $M$, then, for some constant $C>0$ depending only on $(M,g)$,
$$ W_1\left(\phi^+dx, \phi^- dx\right) \le\frac{C}{\sqrt\lambda} \|\phi\|^{}_{L^1(M)}.$$
\end{thm}
This inequality is sharp. We recall the basic intuition that a Laplacian eigenfunction may, at scale $\lambda^{-1/2}$ (the wavelength), be understood as a random wave. This suggests that one has to move mass at least a distance comparable to the wavelength and examples on the torus $\mathbb{T}^d$ or the sphere $\mathbb{S}^d$ show that this is indeed the case.

\subsection{Spectral Bounds on Transport}
The purpose of this section is to show that a variation of this result exists on finite graphs; we will prove this for the Earth Mover's Distance $p=1$.

\begin{theorem} Let $M= I- A D^{-1}$ and let $M\phi_k = \lambda_k \phi_k$. Then $0 \leq \lambda_k \leq 2$ and
$$ W_1( \phi_k^+, \phi_k^-) \leq \frac{1}{1-\left|1-\lambda_k\right|} \|\phi_k\|_{\ell^1}.$$
\end{theorem}
Note that, since $\phi=\phi_k^+-\phi_k^-$ has mean 0, the measures $\phi_k^+$ and $\phi_k^-$ have the same mass, and thus one can be transported to the other.  When we consider the asymptotic behavior of $W_1(\phi_k^+,\phi_k^-)$ on cycle graphs of increasing size, we see that this bound is a natural analog of Peyre's result \cite{pey} to graphs---it scales sharply with respect to $\dot H^{-1}(\phi_k)$ (see \S 3.4). We observe that this bound degenerates if $|\lambda_k - 1|$ is close to 1 and this a consequence of the proof. We also note that we always have the trivial transport inequality
$$ W_1( \phi_k^+, \phi_k^{-}) \leq \mbox{diam}(G) \left\lVert\phi_k^+\right\rVert_{\ell^1}=\frac{\mbox{diam}(G)}{2}\left\lVert\phi_k\right\rVert_{\ell^1},$$ and thus the bound in the Theorem is preferable to the trivial bound only when $$|1-\lambda_k|<1-\frac{2}{\text{diam}(G)}.$$

In many of the interesting cases for applications (graphs with good mixing properties), we can expect a spectral gap that quantitatively bounds $|\lambda_2 - 1| < 1$.

\subsection{Applying the Theorem to obtain Transport Bounds}

For an arbitrary distribution $\mu$, we may use this bound to measure the Wasserstein distance to the uniform distribution on a graph with $n$ vertices.  The observation above motivates splitting $\mu$ into mid-range and extreme-frequency components, $$\overline\mu=\sum_{|1-\lambda_k|<1-2/\text{diam}(G)} \left\langle \mu, \phi_k \right\rangle \phi_k$$
 and $$\underline\mu=\sum_{|1-\lambda_k|\ge 1-2/\text{diam}(G)}  \left\langle \mu, \phi_k \right\rangle \phi_k.$$ We then transport $\overline\mu$ by propagating infinitely, and bound $\underline\mu$ with a diameter bound:\begin{align*}W_1\left(\mu,dx\right)&=W_1(\overline\mu^++\underline\mu^+,\overline\mu^-+\underline\mu^-+dx)\\&\le W_1\left(\overline\mu^+,\overline\mu^-\right)+W_1\left(\underline\mu^+,\underline\mu^-+dx\right)\\&\le \sum_{i=0}^\infty \|A^iD^{-i}\overline\mu\|_{\ell^1}+\frac{\mbox{diam}(G)}{2}\|\underline\mu-dx\|_{\ell^1}\end{align*}
 While the above spectral bound is not guaranteed to be smaller than the diameter bound, empirical evidence suggests that it is in general a stronger bound, particularly for graphs with large diameter.  We can test the quality of this bound by using linear programming \cite{peyre} to compute exact Wasserstein distances on graphs.  In \S4, we display the results of doing so with $$\mu=\frac{1}{k}\sum_{j=1}^k \delta_{x_j}$$ using the algorithm to pick the $x_j$: $$ x_{k+1} = \arg\min_{x\in V} \left(  L^{-1/2} \sum_{j=1}^{k}{\delta_{x_j}}\right)(x),$$ compared against picking vertices $x_k$ uniformly at random (without repetition) and averaging over 1000 Wasserstein distances obtained in this manner.  This approach yields promising computational results across a number of large graphs.

\subsection{A Case Study: Cycle Graphs}
In this subsection, we will look carefully at the behavior of the cycle graphs $C_n$ in the context of the above result.  We recall the eigenvectors of $M=\frac{1}{2}L$ on $C_n$ are  $$\phi_{k}(x)=(n/2)^{-1/2}\cos\left(\frac{2\pi kx}{n}\right)$$ and $$\phi_{n-k}(x)=(n/2)^{-1/2}\sin\left(\frac{2\pi kx}{n}\right),$$ for $0<k<n/2$, with $\phi_0(x)\equiv n^{-1/2}$ and, if $n$ is even, $\phi_{n/2}(x)=n^{-1/2}(-1)^x$ and corresponding eigenvalues $$\lambda_k=1-1\cos\left(\frac{2\pi k}{n}\right).$$
Note that $\phi_0$ is the constant eigenvector here and, since cosine is even, $\lambda_k=\lambda_{n-k}$ for all $k\ne 0$.
We have elected to use the real eigenvectors in order to apply our arguments, though it is worth noting that we can change basis for the dimension 2 eigenspaces and simply write
$$\phi_k(x)=n^{-1/2}\exp\left(\frac{2\pi ikx}{n}\right),$$  where $0\le k<n$ with all $\lambda_k$ the same as above.  Then, $$\frac{1}{1-|1-\lambda_k|}=\frac{1}{1-\cos\left(\frac{2\pi k}{n}\right)}\approx 2\left(\frac{n}{2\pi k}\right)^2$$ for small $k$, by Taylor expansion.  Further, $$\|\phi_k\|_{\ell^1}\approx\|\phi_{n-k}\|_{\ell^1}\approx (n/2)^{-1/2}\int_{0}^n \left|\cos\left(\frac{2\pi kx}{n}\right)\right|dx=\frac{\sqrt{8n}}{\pi}.$$
On the other hand, $W_1(\phi_k^+,\phi_k^-)\approx W_1(\phi_{n-k}^+,\phi_{n-k}^-)$ can be approximated by the continuous analogue, where it is clear from symmetry that the optimal way to transport the sine wave is sending all mass to the nearest zero, where the positive and negative mass will cancel.  This endures a cost of $$(n/2)^{-1/2}4k\int_0^{n/4k}x\sin\left(\frac{2\pi kx}{n}\right)dx=(n/2)^{-1/2}4k\left(\frac{n}{2\pi k}\right)^2,$$ integrating by parts.  Putting it all together, this yields $$W_1(\phi_k^+,\phi_k^-)\approx \frac{\pi k}{n}\cdot\frac{1}{1-|1-\lambda_k|}\|\phi_k\|_{\ell^1}.$$  We may let $k\le n/100$ so that $k$ is small enough for the Taylor expansion to be good, but nonetheless on the order of $n$: then we see the bound in Theorem 2 is sharp up to constants.  Observe that when we apply the fractional inverse Laplacian $ L^{-1/2}$ to the point mass $\delta_0$, we get 
\begin{align*} L^{-1/2}(\delta_0)=\sum_{k=1}^{n-1} \frac{\phi_k(0)}{\lambda_k^{1/2}}\phi_k&=\sum_{k=1}^{n-1} \frac{1}{\left(2-2\cos\left(\frac{2\pi k}{n}\right)\right)^{1/2} n}\exp\left(\frac{2\pi ikx}{n}\right)
\\&=2\sum_{k=1}^{\lfloor n/2\rfloor}\frac{1}{\left(2-2\cos\left(\frac{2\pi k}{n}\right)\right)^{1/2} n}\cos\left(\frac{2\pi kx}{n}\right)
\\&\approx \frac{1}{\pi}\sum_{k=1}^{\lfloor n/2\rfloor} \frac{1}{k}\cos\left(\frac{2\pi k x}{n}\right),\end{align*} for all odd $n$, again approximating with Taylor expansion.  (If $n$ is even, we will get an extra $(-1)^x/n\sqrt{2}$ term corresponding to $\phi_{n/2}$, but this will vanish in the limit we are about to take.)  We caution the reader that the $\lambda_k$ above are the eigenvalues of $L=2M$, and are thus double the eigenvalues of $M$ referred to in the preceding computations.  Rescaling with $x=n\theta/2\pi$, we have $$ L^{-1/2}(\delta_0)\approx \frac{1}{\pi}\sum_{k=1}^{\lfloor n/2\rfloor} \frac{1}{k}\cos\left(k\theta\right).$$  Fixing $\theta$ and letting $n\to\infty,$ this is simply the Fourier series for $$\lim_{n\to\infty}L^{-1/2}(\delta_0)=-\frac{1}{\pi}\ln \left|2\sin\left(\theta/2\right)\right|.$$  Taking this limit and rescaling really is just transitioning us to the continuous setting: the eigenfunctions of the Laplacian on the unit circle $\mathbb S^1$ identified with $[0,2\pi)/\sim$ are $(2\pi)^{-1/2}\exp(ikx)$, with eigenvalue $k^2$, for $k\in\mathbb{Z}$.
So the fractional inverse Laplacian $ L^{-1/2}$ of a point mass on the circle is
\begin{align*} L^{-1/2}(\delta_0)&=\sum_{k\ne 0} \frac{\phi_k(0)}{\lambda_k^{1/2}}\phi_k=\sum_{k\ne 0} \frac{1}{2\pi|k|}\exp(ikx)
\\&= \frac{1}{\pi}\sum_{k=1}^\infty \frac{1}{k}\cos(kx)
=-\frac{1}{\pi}\ln|2\sin(x/2)|,\end{align*} precisely our function in the discrete case.

\begin{figure}[!htb]\centering
\includegraphics[scale=.7]{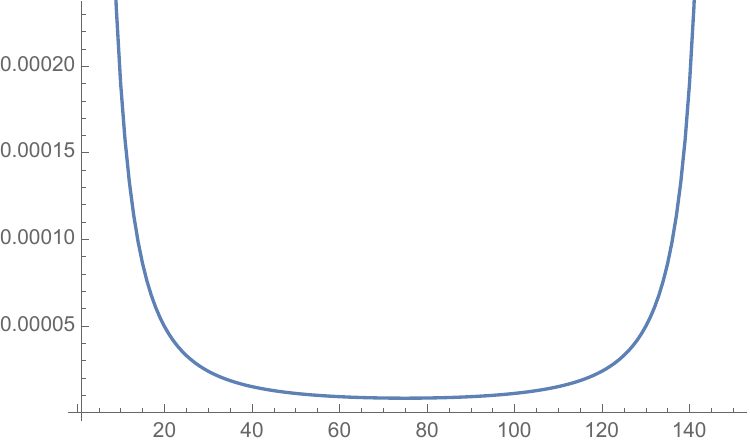}
\caption{The difference between $ L^{-1/2}(\delta_0)$ on $\mathbb S^1$ and the fractional Laplacian $ L^{-1/2}(\delta_0)$ on $C_{150}$ is small.} 
\end{figure}

In Figure 6 we show the difference between the inverse fractional Laplacian on a point mass $ L^{1/2}(\delta_0)$ in the continuous and discrete settings.  The two outputs are almost identical and thus their difference is quite small (this holds even for $very$ small values of $n$).  It is worth recalling the recent Theorem of Pausinger, which proves that this algorithm belongs to a large class which produces the van der Corput sequence on the torus (and thus achieves optimal discrepancy, up to constants) \cite{pausinger}.  Further, Steinerberger's recent result \cite{steinii} indicates that this algorithm performs extremely well on the torus, and by extension large cycle graphs, and that the sequence of points rapidly becomes very evenly distributed.

\subsection{A Case Study: Torus Grid Graphs}
In this subsection, we examine another class of graphs: torus grid graphs.  The $m\times n$ torus grid graph $T_{m,n}$ is the Cartesian product of $C_m$ and $C_n$, so we can apply many of our computations from the previous subsection here.  In particular, the eigenvectors for $T_{m,n}$ are the Kronecker products of pairs of eigenvectors $(\phi_j,\phi_k)$ from $C_m$ and $C_n$, respectively, with corresponding eigenvalues under $M=\frac{1}{4}L$ of \begin{align*}\lambda_{j,k}&=\frac{\lambda_j+\lambda_k}{2}=1-\frac{1}{2}\left(\cos\left(\frac{2\pi j}{m}\right)+\cos\left(\frac{2\pi k}{n}\right)\right).\end{align*}  (The factor of $1/2$ appears because cycle graphs are 2-regular while torus grid graphs are 4-regular, and vanishes if we instead use $L$ here.)  Then we have
\begin{align*}\frac{1}{1-|1-\lambda_{j,k}|}&=\frac{1}{1-\frac{1}{2}\left(\cos\left(\frac{2\pi j}{m}\right)+\cos\left(\frac{2\pi k}{n}\right)\right)}\approx\frac{4}{\left(\frac{2\pi j}{m}\right)^2+\left(\frac{2\pi k}{n}\right)^2}\end{align*} for small $j,k$, by Taylor expansion.  Further, $$\|\phi_{j,k}\|_{\ell^1}=\|\phi_j\|_{\ell^1}\|\phi_k\|_{\ell^1}\approx \frac{8\sqrt{nm}}{\pi^2}.$$
To find $W_1(\phi_{j,k}^+,\phi_{j,k}^-)$, note that the support of $\phi_{j,k}^+$ consists of checkerboarded rectangles, and the most efficient way to transport the positive mass to the negative mass will be along the higher frequency direction---that is, horizontally if $m/j<n/k$, and vertically otherwise.  Without loss of generality, let us suppose we are in the former case.  Then this incurs a cost of \begin{align*}\sum_{i=0}^{n-1}\lvert\phi_k(i)\rvert W_1(\phi_j^+,\phi_j^-)&=\|\phi_k\|_{\ell^1}W_1(\phi_j^+,\phi_j^-)\approx \frac{\sqrt{8n}}{\pi}\cdot(m/2)^{-1/2}4j\left(\frac{m}{2\pi j}\right)^2.\end{align*}  Applying our assumption that $m/j<n/k$ to our earlier estimate, we see \begin{align*}\frac{1}{1-|1-\lambda_{j,k}|}&\approx \frac{4}{\left(\frac{2\pi j}{m}\right)^2+\left(\frac{2\pi k}{n}\right)^2}\le \frac{4}{\left(\frac{2\pi k}{n}\right)^2+\left(\frac{2\pi k}{n}\right)^2}=2\left(\frac{n}{2\pi k}\right)^2.\end{align*}  Putting it all together, \begin{align*}\frac{\sqrt{8n}}{\pi}(m/2)^{-1/2}4j\left(\frac{m}{2\pi j}\right)^2&\approx W_1(\phi_{j,k}^+,\phi_{j,k}^-)\\&\le\frac{1}{1-|1-\lambda_{j,k}|}\|\phi_{j,k}\|_{\ell^1}\le 2\left(\frac{n}{2\pi k}\right)^2 \cdot \frac{8\sqrt{nm}}{\pi^2},\end{align*} and thus, taking the quotient of the two sides in the above inequality our bound is off by (at most) a factor of $\pi k^2 m(n^{2}j)^{-1}.$  Note that when $k/n=j/m$ (i.e., when the horizontal and vertical components of $\phi_{j,k}$ have the same frequency), this simplifies to $\pi k/n$, precisely our result on cycle graphs.

\section{Numerics}
Below we provide numerics on a variety of graphs demonstrating the performance of the algorithm and the bound from Theorem 2.  In particular, we compare the performance of vertices selected according to our algorithm against that of randomly selected vertices.  The differences between our vertex sequences and random vertex sequences may seem marginal, but this is partly due to the fact that the diameter of some of these graphs is quite small. For instance, the Truncated Tetrahedral graph, with diameter 3, only has 12 vertices, so we will hardly be able to distinguish the performance of the algorithm's vertices from randomly selected vertices on such a small set---it is impressive that we see a difference at all.  We see that for the Faulkner-Younger Graph and the Level 2 Menger Sponge the difference becomes significantly more drastic.  Many of the graphs are quite well connected, which makes the transport problem easier than on sparse graphs (e.g., on complete graphs it makes no difference at all which vertices are selected, the transport cost only depends on the number of vertices).  In all the tables in this section, $x_j$ were computed directly using the recursive definition for the algorithm ($\alpha=.5$) given in Section 2 with any ties broken randomly, and the exact Wasserstein distances $$W_1\left(\frac{1}{k}\sum_{j=1}^k \delta_{x_j},dx\right)$$ were subsequently computed using the dual linear program in \cite{peyre} in the ``Algorithm" row.  For graphs which are not vertex-transitive, the performance of the algorithm depends upon the arbitrary initial vertex chosen, and thus all choices of initial vertex were attempted and the transport costs averaged.  In the ``Random" row, 1000 uniformly randomly selected sets of $k$ distinct vertices were taken, and the corresponding Wasserstein distances were averaged.

\begin{figure}[H]\centering
\begin{minipage}{.4\textwidth}\centering
\includegraphics[width=\linewidth]{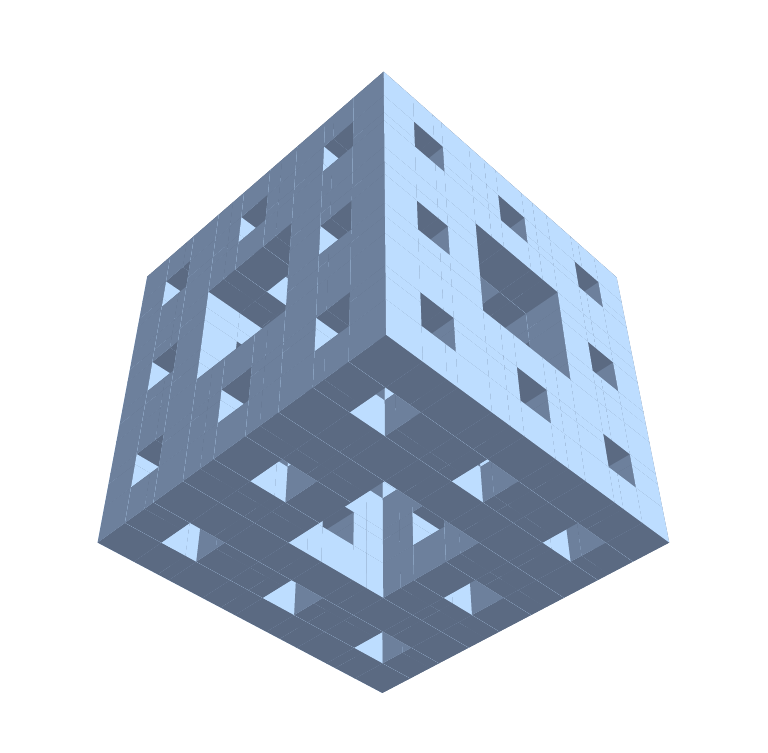}
\captionof{figure}{The Menger Sponge, whose 400 cubes form the vertices of a connectivity graph.} \end{minipage}%
\begin{minipage}{.6\textwidth}\centering
\includegraphics[width=\linewidth]{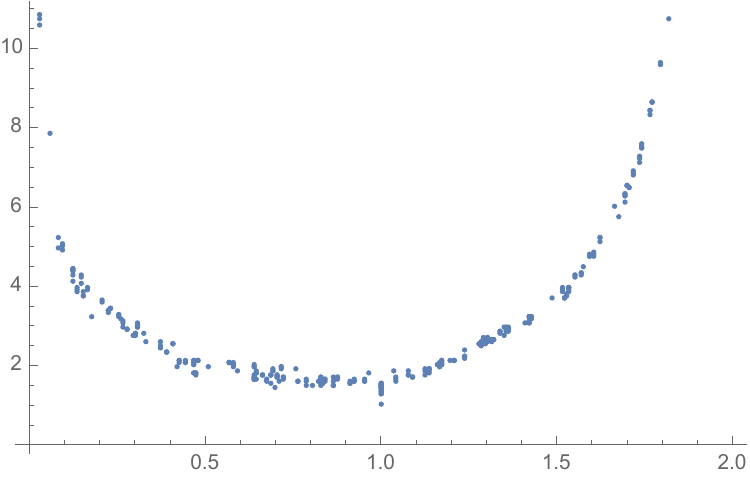}
\captionof{figure}{The tightness of the bound in Theorem 2, applied to the eigenfunctions of the Level 2 Menger Sponge Connectivity Graph (computed as a quotient of the right and left sides).} 
\end{minipage}
\end{figure}
 
\subsection{Connectivity Graph of Level 2 Menger Sponge}
The Level 2 Menger sponge is the object obtained beginning with a cube and drilling out the middle square of each face (viewed as a three by three grid of squares), and then iterating this process one more time on the smaller cubes (see Fig. 7).  We can then generate a connectivity graph of the remaining 400 smaller cubes (each one ninth the side length of the original cube).  Note that this is not a regular graph.  In Figure 8, we see the that, on the Level 2 Menger Sponge Connectivity Graph, the Theorem 2 bound is tightest for mid-range eigenvalues.  This is to be expected, due to the blow-up of the $1/(1-|1-\lambda|)$ term at the extremes, where a diameter bound is tighter (see $\S 3.2$).  But we see here that, even for very small eigenvalues, the bound is fairly tight.

\begin{table}[H]\caption{$W_1(\mu,dx)$ for the Connectivity Graph of a Level 2 Menger Sponge}\centering
\begin{tabular}{c|cccccccc}
\text{No. of vertices}& 1   & 3    & 5    & 10   & 15   & 20   & 25   & 30   \\
\text{Algorithm} & 9.94 & 6.48 & 5.01 & 3.54 & 2.92 & 2.55 & 2.31 & 2.11 \\
\text{Random} &9.94  & 6.73 & 5.52 & 4.25 & 3.63 & 3.20 & 2.90 & 2.69
\end{tabular}
\end{table}

\begin{figure}[H]\centering\begin{minipage}{.5\textwidth}\centering
\begin{tikzpicture}[scale=.8]
\foreach \a in {1,2,...,9}{
\draw [thick] (\a*360/9: 3cm) --  (\a*360/9 + 360/9: 3cm);
\filldraw [thick] (\a*40:3cm) circle (0.08cm);
};
\filldraw [thick] (0.8, -0.4) circle (0.08cm);
\filldraw [thick] (-0.8, -0.4) circle (0.08cm);
\filldraw [thick] (0, 1) circle (0.08cm);
\draw [thick] (0.8, -0.4) --  (8*360/9: 3cm);
\draw [thick] (0.8, -0.4) -- (-0.8, -0.4) -- (0,1) -- (0.8, -0.4);
\draw [thick] (-0.8, -0.4) --  (5*360/9: 3cm);
\draw [thick] (0,1) --  (2*360/9: 3cm);
\draw [thick] (1*360/9: 3cm) --  (3*360/9: 3cm);
\draw [thick] (4*360/9: 3cm) --  (6*360/9: 3cm);
\draw [thick] (7*360/9: 3cm) --  (9*360/9: 3cm);
\node at (80:3.35cm) {\Large 1};
\node at (-80:3.35cm) {\Large 2};
\node at (-.8,-.75) {\Large 3};
\node at (160:3.35cm) {\Large 4};
\node at (0:3.35cm) {\Large 5};
\node at (.8,-.75) {\Large 6};
\node at (120:3.35cm) {\Large 7};
\node at (-120:3.35cm) {\Large 8};
\node at (-.2,1.2) {\Large 9};
\node at (-40:3.35cm) {\Large {10}};
\node at (40:3.35cm) {\Large {11}};
\node at (-160:3.35cm) {\Large {12}};
\end{tikzpicture}
\caption{The sequence of vertices picked by the algorithm on the Truncated Tetrahedral Graph.} 
\end{minipage}%
\begin{minipage}{.5\textwidth}\centering
  \begin{tikzpicture}[scale=.8]
\foreach \a in {1,2,...,7}{
\filldraw (\a*360/7: 3cm) circle (0.08cm);
};
\foreach \a in {1,2,...,7}{
\draw [thick] (\a*360/7: 3cm) --  (\a*360/7 + 360/7: 3cm);
};
\filldraw [thick] (0.8, 0) circle (0.08cm);
\draw [thick] (0.8, 0) --  (1*360/7: 3cm);
\draw [thick] (0.8, 0) --  (7*360/7: 3cm);
\draw [thick] (0.8, 0) --  (0.2, 1.5);
\filldraw [thick] (0.2, 1.5) circle (0.08cm);
\filldraw [thick] (0,0) circle (0.08cm);
\filldraw [ thick] (-1,0) circle (0.08cm);
\filldraw [ thick] (0,-1.5) circle (0.08cm);
\draw [ thick] (0.2, 1.5) --  (2*360/7: 3cm);
\draw [ thick] (0.2, 1.5) --  (0,0);
\draw [ thick] (-1, 0) --  (0,0);
\draw [ thick] (0, -1.5) --  (0,0);
\draw [ thick] (-1, 0) --  (3*360/7: 3cm);
\draw [ thick] (-1, 0) --  (4*360/7: 3cm);
\draw [ thick] (0,-1.5) --  (5*360/7: 3cm);
\draw [ thick] (0,-1.5) --  (6*360/7: 3cm);
\node at (1*360/7: 3.35cm) {\Large 1};
\node at (4*360/7: 3.35cm) {\Large 2};
\node at (-0.3,-1.4) {\Large 3};
\node at (0.5,1.5) {\Large 4};
\node at (-1.4,0) {\Large 5};
\node at (7*360/7: 3.35cm) {\Large 6};
\node at (5*360/7: 3.35cm) {\Large 7};
\node at (2*360/7: 3.35cm) {\Large 8};
\node at (-0.2,-0.3) {\Large 9};
\node at (0.7,-0.3) {\Large 10};
\node at (3*360/7: 3.35cm) {\Large 11};
\node at (6*360/7: 3.35cm) {\Large 12};
   \end{tikzpicture}
\caption{The sequence of vertices picked by the algorithm on the Frucht Graph.} 
\end{minipage}
\end{figure}

\subsection{Truncated Tetrahedral Graph}
The Truncated Tetrahedral Graph (see Fig. 9) is a 3-regular, vertex-transitive graph on 12 vertices.  It is the 1-skeleton of the Archimedean solid formed by truncating each vertex of a tetrahedron.
\begin{table}[H]\caption{$W_1(\mu,dx)$ for the Truncated Tetrahedral Graph}\centering
\begin{tabular}{c|cccccccccc}
\text{No. of vertices} &1    & 2    & 3    & 4    & 5    & 6    & 7    & 8    & 9    & 10   \\
\text{Algorithm} & 1.92 & 1.17 & 0.83 & 0.67 & 0.58 & 0.50 & 0.42 & 0.33 & 0.28 & 0.23 \\
 \text{Random} &1.92 & 1.35 & 1.01 & 0.84 & 0.72 & 0.58 & 0.52 & 0.43 & 0.34 & 0.27
\end{tabular}
\end{table}

\subsection{Frucht Graph}
The Frucht Graph (see Fig. 10) is a 3-regular graph on 12 vertices, and has trivial automorphism group despite being degree-regular.

\begin{table}[H]\caption{$W_1(\mu,dx)$ for the Frucht Graph}\centering
\begin{tabular}{c|cccccccccc}
\text{No. of vertices} & 1    & 2    & 3    & 4    & 5    & 6    & 7    & 8    & 9    & 10   \\
\text{Algorithm} & 1.93 & 1.17 & 0.86 & 0.67 & 0.59 & 0.50 & 0.42 & 0.34 & 0.29 & 0.23  \\
\text{Random} & 1.93 & 1.34 & 1.04 & 0.85 & 0.73 & 0.59 & 0.52 & 0.43 & 0.35 & 0.27
\end{tabular}
\end{table}

\subsection{Faulkner-Younger Graph}
The Faulkner-Younger Graph on 44 vertices (see Fig. 11) is a 3-regular non-Hamiltonian graph (that is, there is no path along its edges that traverses every vertex exactly once).
\begin{figure}[H]\centering
\includegraphics[scale=.7]{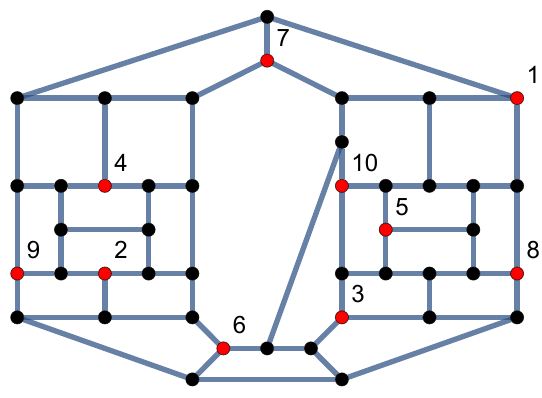}
\caption{The first ten vertices picked by the algorithm on the Faulkner-Younger Graph.  Each label is above and to the right of the corresponding vertex.} 
\end{figure}

\begin{table}[H]\caption{$W_1(\mu,dx)$ for the Faulkner-Younger Graph}\centering
\begin{tabular}{c|cccccccccc}
\text{No. of vertices} & 1  & 2    & 3    & 4    & 5    & 6    & 7    & 8    & 9    & 10   \\
\text{Algorithm} & 4.17 & 2.67 & 2.05 & 1.71 & 1.52 & 1.34 & 1.23 & 1.15 & 1.05 & 0.97 \\
\text{Random} & 4.17 & 3.08 & 2.57 & 2.24 & 2.01 & 1.83 & 1.68 & 1.56 & 1.46 & 1.37
\end{tabular}
\end{table}

\subsection{Erd\H os-R\'enyi Random Graphs}
The Erd\H os-R\'enyi model for random graphs $G(n,p)$ is given by including an edge between each pair of the $n$ vertices independently with probability $p$.  Here we display the performance of the algorithm on two such graphs, one taken from $G(100,.06)$ (a sparse graph, see Fig. 12) and another from $G(100,.2)$ (a dense graph, see Fig. 13).

\begin{table}[H]\caption{$W_1(\mu,dx)$ for Erd\H os-R\'enyi Random Graphs}\centering
\begin{tabular}{c|cccccccccc}
\text{No. of vertices} & 1  & 3  & 5 & 10  & 15  & 20 & 25 & 30  \\
\text{Algorithm, Sparse Graph} & 2.72 & 2.61 & 2.13 & 1.52 & 1.22 & 1.02 & 0.85 & 0.74 \\
\text{Random, Sparse Graph} & 2.72 & 2.11 & 1.82 & 1.44 & 1.22 & 1.06 & 0.93 & 0.83 \\
\text{Algorithm, Dense Graph} &1.79 & 1.53 & 1.31 & 1.01 & 0.86 & 0.80 & 0.75 & 0.70 \\
\text{Random, Dense Graph} & 1.79 & 1.46 & 1.26 & 0.99 & 0.88 & 0.81 & 0.75 & 0.70
\end{tabular}
\end{table}

It is no surprise that the dense graph exhibits little variation between the transport cost of random vertices and of the algorithm's---after all, any pair of vertices has many short paths between them.  In fact, this particular graph has diameter 3.  Thus, for sufficiently dense graphs it is largely irrelevant which vertices are selected: the transport cost will be low.  The sparse graph displayed has diameter 6 and is thus more interesting: while random vertices initially outperform the algorithm, by 15 vertices selected they are matched, after which the algorithm surpasses the random vertices.  That is, even in highly irregular graphs such as this one where random vertices perform well at first, the algorithm nonetheless manages to catch up even with a relatively small number of vertices.

\begin{figure}[h!]\centering
\begin{minipage}{.5\textwidth}\centering
\includegraphics[width=.9\linewidth]{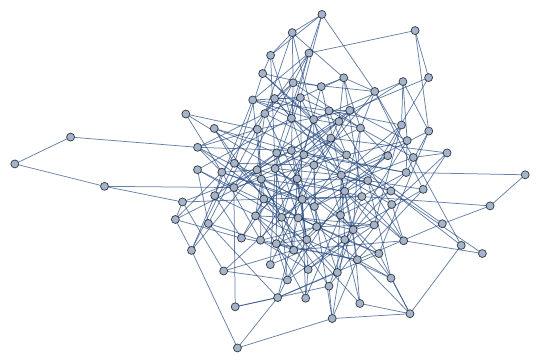}
\captionof{figure}{A sparse Erd\H os-R\'enyi Graph from $G(100,.06)$.} \end{minipage}%
\begin{minipage}{.5\textwidth}\centering
\includegraphics[scale=.5]{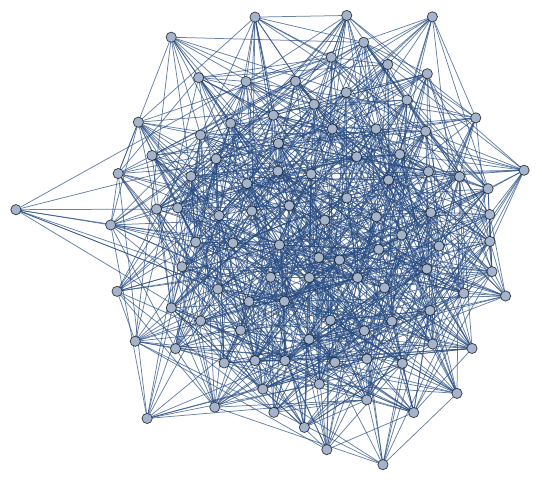}
\captionof{figure}{A dense Erd\H os-R\'enyi Graph from $G(100,.2)$.} 
\end{minipage}
\end{figure}
 
\subsection{Complete 3-ary Tree}
The complete 3-ary tree of depth 4 is a rooted tree, where each vertex has 3 children, except for the fourth generation of vertices which all have no children, yielding a total of 40 vertices (see Fig. 14).
\begin{figure}[H]\centering
\includegraphics[scale=1]{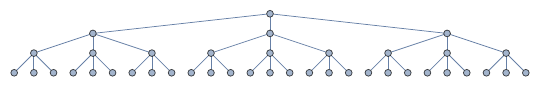}
\caption{The 3-ary tree of depth 4} 
\end{figure}

\begin{table}[H]\caption{$W_1(\mu,dx)$ for the 3-ary tree of depth 4}\centering
\begin{tabular}{c|cccccccccc}
\text{No. of vertices} & 1 & 2 & 3 & 4 & 5 & 6 & 7 & 8 & 9 & 10  \\
\text{Algorithm} &4.25 & 3.58 & 2.75 & 2.83 & 2.55 & 2.07 & 1.98 & 1.73 & 1.34 & 1.37 \\
\text{Random} & 4.25 & 3.62 & 3.12 & 2.93 & 2.74 & 2.48 & 2.33 & 2.17 & 1.99 & 1.86 \\
\end{tabular}
\end{table}

\section{Proofs}
\subsection{Proof of Theorem 1}
\begin{proof}
Note first that, since, for all $v\in\mathbb R^n$, $L^{-2\alpha}v$ has mean 0 (being spanned by $\phi_i$, $i>1$, and thus orthogonal to the constant $\phi_1$), we will always have 
$$\min_{x\in V} \left(  L^{-2\alpha} \sum_{j=1}^{k}{\delta_{x_j}}\right)(x)<0.$$
Using the $\ell^2$ norm $$\|v\|_{\ell^2}^2=\sum_{i=1}^n v_i^2,$$ we observe \begin{align*}\left\lVert L^{-\alpha}(k\mu_k)\right\rVert^2_{\ell^2}=&\left\lVert L^{-\alpha}\left((k-1)\mu_{k-1}\right)\right\rVert^2_{\ell^2}+\left\lVert L^{-\alpha}\left(\delta_{x_k}\right)\right\rVert^2_{\ell^2}\\&+2\left\langle  L^{-\alpha}((k-1)\mu_{k-1}), L^{-\alpha}\left(\delta_{x_k}\right)\right\rangle
\\=&\left\lVert L^{-\alpha}\left((k-1)\mu_{k-1}\right)\right\rVert^2_{\ell^2}+\left\lVert L^{-\alpha}\left(\delta_{x_k}\right)\right\rVert^2_{\ell^2}\\&+2\left\langle  L^{-2\alpha}((k-1)\mu_{k-1}),\delta_{x_k}\right\rangle,\end{align*}
since $ L^{-\alpha}$ is self-adjoint.  Rewriting the inner product term, $$\left\langle  L^{-2\alpha}((k-1)\mu_{k-1}),\delta_{x_k}\right\rangle=\left( L^{-2\alpha}\sum_{j=1}^{k-1}\delta_{x_j}\right)(x_k).$$  But $x_k$ was chosen by the algorithm specifically to minimize that quantity---thus, it is certainly less than the average value of 0, and so $$\left\lVert L^{-\alpha}(k\mu_k)\right\rVert^2_{\ell^2}\le\left\lVert L^{-\alpha}\left((k-1)\mu_{k-1}\right)\right\rVert^2_{\ell^2}+\left\lVert L^{-\alpha}\left(\delta_{x_k}\right)\right\rVert^2_{\ell^2}.$$
Then, by induction, we obtain the desired inequality:
$$\sum_{i=2}^n\frac{|\langle\mu_k,\phi_i\rangle|^2}{\lambda_i^\alpha}=\left\lVert L^{-\alpha}(\mu_k)\right\rVert_{\ell^2}^2\le \left(\max_{j\le k}\left\lVert L^{-\alpha}\left(\delta_{x_j}\right)\right\rVert_{\ell^2}^2\right)k^{-1}.$$

\end{proof}
\subsection{Proof of Theorem 2}
\begin{proof} We recall that $AD^{-1}$ can be interpreted as the propagator of the random walk on the graph $G=(V,E)$. Moreover, we have 
$$ AD^{-1}\phi_k = (1-\lambda_k) \phi_k$$
and observe that $|1 - \lambda_k| \leq 1$.  We proceed in a similar manner to \cite{stein} and interpret diffusion on the graph as one of many ways to transport mass. In particular, we will apply the elementary estimate
$$ W_1(A D^{-1} v, v) \leq \|v\|_{\ell^1}=\sum_{j=1}^n |v_j|$$ to $v=((AD^{-1})^i\phi_k)^\pm$,
$$ W_1(AD^{-1}((AD^{-1})^i\phi_k)^\pm), ((AD^{-1})^i\phi_k)^\pm)\leq \|((AD^{-1})^i\phi_k)^\pm\|_{\ell^1}=\frac{\left|1-\lambda_k\right|^i}{2}\|\phi_k\|_{\ell^1}.$$
In particular, we will transport $\phi_k$ to $(AD^{-1})^m \phi_k$ through its positive and negative parts after each diffusion.  For large $m$, this measure almost vanishes since
$$(AD^{-1})^m\phi_k =\left(1-\lambda_k\right)^{m} \phi_k.$$
We perform this operation until some arbitrary $m$ and then use the trivial bound on the remaining measure.
This shows that the total transport can be bounded by
\begin{align*} W_1( \phi_k^+, \phi_k^-) &\leq    \frac{\left| 1 - \lambda_k\right|^{m} \mbox{diam}(G)}{2} \|\phi_k\|_{\ell^1} + \sum_{i=0}^{m-1} \left|1-\lambda_k\right|^i\|\phi_k\|_{\ell^1}\\&=\left( \frac{\left| 1 - \lambda_k\right|^{m} \mbox{diam}(G)}{2} + \frac{1-\left|1-\lambda_k\right|^m}{1-\left|1-\lambda_k\right|}\right)\|\phi_k\|_{\ell^1}\\&=\left(\left|1-\lambda_k\right|^m\left[\frac{\mbox{diam}(G)}{2}-\frac{1}{1-\left|1-\lambda_k\right|}\right]+\frac{1}{1-\left|1-\lambda_k\right|}\right)\|\phi_k\|_{\ell^1} .\end{align*}
We observe that this bound is monotonic in $m$, with direction depending on the sign of the bracketed expression.  If $$|1-\lambda_k|\ge 1-\frac{2}{\mbox{diam}(G)},$$ the bound is monotonically increasing and we set $m=0$, recovering the initial diameter bound.  If $$|1-\lambda_k|<1-\frac{2}{\mbox{diam}(G)},$$ we have a monotonically decreasing bound, and take the limit as $m\to\infty$, yielding the desired bound of $$\frac{1}{1-|1-\lambda_k|}\|\phi_k\|_{\ell^1}.$$
\end{proof}

\section{Connection to other Results}
\subsection{Low-discrepancy point sets.}
A classical problem in the study of \textit{irregularities of distribution} is to construct sequences $(x_n)_{n=1}^{\infty}$ on the unit interval $[0,1]$ such that $\left\{x_1, \dots, x_n\right\}$ is fairly evenly distributed over the unit interval for all $n \in \mathbb{N}$. The problem has now been solved completely: it is known, this is a result of Schmidt \cite{schmidt}, that for any sequences on $[0,1]$, there exists infinitely many $n \in \mathbb{N}$ such that
$$  \max_{0 \leq x \leq 1} \left| \frac{\# \left\{1 \leq i \leq n: x_i \leq x\right\} }{n} - x \right| \geq \frac{1}{100}\frac{\log{n}}{n}.$$
The quantity on the left-hand side is also known as discrepancy (or extreme discrepancy, $L^{\infty}-$discrepancy), we refer to the textbook of Dick \& Pillichshammer \cite{dick}.
Steinerberger recently established \cite{steini,steinii} that a greedy sequence defined via
$$ x_{n+1} = \arg\min_{x\in [0,1]} \left(  L^{-1/2} \sum_{k=1}^{n}{\delta_{x_k}}\right)(x)$$
satisfies excellent distribution properties for all $n \in \mathbb{N}$---numerical examples show that the arising sequences seem to be remarkably close to the best possible bound $n^{-1}\log{n}$ (down to the
level of the constant). The argument is somewhat different and uses the Koksma-Hlawka inequality \cite{koksma} and classical Fourier Analysis. In particular, this result is stronger than what is guaranteed by Theorem 1.

\subsection{Greedy minimization}
Greedy minimizations such as the one explored in this paper are well-behaved in general.  In \cite{brown2}, Steinerberger and the author showed that for any $f:\mathbb{T}\to\mathbb{R}$ with $\widehat f(k)\ge c|k|^{-2}$, the sequence defined by $$x_n=\arg\min_{x\in \mathbb T} \sum_{k=1}^{n-1} f(x-x_k)$$ is well-distributed.  Note that this algorithm is summing shifted copies of $f$ and finding the smallest value, essentially ``filling in the gaps" in the point set.  Independently of initial choice of $\{x_1,\dots,x_k\}$, we have $$W_2\left(\frac{1}{n}\sum_{i=1}^n \delta_{x_i},dx\right)\lesssim \frac{c}{\sqrt{n}}.$$
\\In higher dimensions, the picture is even nicer: defining an appropriate analogue of the $\widehat f(k)\ge c|k|^{-2}$ condition, we find ourselves minimizing the energy of Green's function-like kernels, and now have, for $d\ge 3$, $$W_2\left(\frac{1}{n}\sum_{i=1}^n \delta_{x_i},dx\right)\lesssim \frac{1}{n^{1/d}}.$$  The exception is dimension $d=2$, where the bound weakens to $\sqrt{\log n}/n^{1/2}$.  It is unknown if this bound is sharp.

\subsection{Leja points} Leja points can be defined, in the utmost level of generality, for any symmetric kernel $k: X \times X \rightarrow \mathbb{R} \cup \left\{ \infty\right\}$ on a compact Hausdorff space. We remain on smooth compact manifolds $M$, a natural example for 
the kernel is 
$$ k(x,y) = \frac{1}{d_g(x,y)^s} \qquad \mbox{where}~d_g(x,y)~\mbox{is the geodesic distance and}~s>0.$$
We can then define, in an iterative fashion, for a given initial point $x_1 \in M$, a sequence $(x_k)_{k=1}^{\infty},$ in such a way that
$$ \sum_{ k=1}^{n-1}{ k(x_n, x_k)} = \inf_{x \in M}  \sum_{ k=1}^{n-1}{ k(x, x_k)}.$$
Put differently, we add a new point $x_n$ in a greedy fashion in such a way that the total energy 
$$ \sum_{ k, \ell =1 \atop k \neq \ell}^{n}{ k(x_k, x_{\ell})} \qquad \mbox{is as small as possible.}$$
These sets were introduced by Edrei \cite{ed} and intensively studied by Leja \cite{leja} after whom they are named. The most commonly used kernel is $k(x,y) = -\log{|x-y|}$ (such that minimizing the sum is the same as maximizing the product of the distances). 
Leja points have a number
of applications in numerical analysis \cite{num0, num1, num15, num2, num3}. Pausinger \cite{pausinger} recently gave a very precise description of Leja sequences on $\mathbb{T}$ for fairly general kernel functions and established a connection to binary digit expansion. For the Riesz kernel $k(x,y) = |x-y|^{-s}$, it is known that Leja sequences are asymptotically uniformly distributed \cite{lopez}.  We are not aware of any study of Leja \textit{vertices} on graphs; while one could take existing kernels, for example $k(x,y) = |x-y|^{-s}$, and consider them on graphs, there is little reason to assume that such vertices will have many special properties: Graphs are simply too flexible.  
 We can summarize the approach in this paper as stating that
\begin{quote}
there is a very good reason to believe (see the Figures in this paper) that considering $k(x,y)$ to be the Green's function of the inverse Laplacian leads to well-distributed sets of vertices.
\end{quote}
Moreover, we are able to analyze the continuous limit of manifolds and are able to obtain a quantitative bound showing that the bounds are more regularly distributed than simply exhibiting uniform distribution.

\subsection{Riesz points} Riesz points refer, at great level of generality, to point sets minimizing energy expressions of the following form
$$ \arg\min_{x\in M^n}\sum_{k, \ell =1 \atop k \neq \ell}^{n}{\frac{1}{\|x_k -x_{\ell}\|^s}}.$$
The problem was first stated on $\mathbb{S}^2$ with $s=1$ by Thomson \cite{thomson} in 1904 and has since inspired a lot of work, we refer to \cite{brauchart, dahlberg, hardin, saff} and references therein. We make a connection with two contributions in particular. The first is due to Beltran, Corral and Criado del Rey \cite{carlos}: they show that if we consider sets of $n$ points on a compact manifold chosen to attain minimal energy
$$\min_{x\in M^n}\sum_{k, \ell =1 \atop k \neq \ell}^{n}{G(x_k, x_{\ell})},$$
where $G$ is the Green's function of the Laplacian on $M$, then the sequence of measures converges weakly to the uniform measure
$$ \frac{1}{n} \sum_{k=1}^{n}{\delta_{x_k}} \rightharpoonup dx.$$
This can be considered the static analogue (since one finds the minimal arrangement for all $n$ points) of our problem (keeping the previous $n-1$ points fixed and then adding the point so as to minimize energy).  The second contribution that we highlight is very recent and due to Marzo \& Mas \cite{marzo}. They studied the specific problem of minimizing the $s-$Riesz energy
$$ E_s =  \sum_{k, \ell =1 \atop k \neq \ell}^{n}{\frac{1}{\|x_k -x_{\ell}\|^s}} \qquad \mbox{on}~\mathbb{S}^d$$
and estimating the spherical cap discrepancy of the minimizing point set: in short, if the points are uniformly distributed, then we would expect the number of points in each spherical cap to be proportional to the volume of the cap; the largest discrepancy is known as spherical cap discrepancy. They use ideas dating back to Wolff: the Riesz energy $E_s$ is comparable to a negative Sobolev norm and, more precisely, there exists a constant $C_{s,d}>1$ such that for all $f \in L^2(\mathbb{S}^d)$,
$$ C_{s,d}^{-1}\|f\|_{H^{(s-d)/2}}^2 \le \int_{\mathbb{S}^d \times \mathbb{S}^d}{ \frac{f(x) f(y)}{\|x-y\|^s} dx dy} \le C_{s,d}\|f\|^2_{H^{(s-d)/2}}.$$
This, while not directly related to our approach, is at least philosophically related: we estimate the Wasserstein distance in negative Sobolev spaces and use the underlying $L^2-$structure.


\begin{thebibliography}{10}

\bibitem{carlos} C. Beltran, N. Corral and J. Criado del Rey, Discrete and Continuous Green Energy on Compact Manifolds, Journal of Approximation Theory 237, 160--185 (2019).

\bibitem{num0} L. Bialas-Ciez, J.-P. Calvi, Pseudo Leja sequences, Annali di Matematica 191, 53--75 (2012).

\bibitem{brauchart} J. Brauchart, Optimal logarithmic energy points on the unit sphere, Math. Comp. 77, 1599-1613 (2008).

\bibitem{brown1} L. Brown, S. Steinerberger, On the Wasserstein Distance between Classical Sequences and the Lebesgue Measure, Trans. Amer. Math. Soc. 373, 8943--8962 (2020).

\bibitem{brown2} L. Brown, S. Steinerberger, Positive-definite Functions, Exponential Sums and the Greedy Algorithm: a Curious Phenomenon, Journal of Complexity (2020).

\bibitem{num1} D. Calvetti, L. Reichel, D. Sorensen, An implicitly restarted Lanczos method for large symmetric eigenvalue problems, Electronic Transactions on Numerical Analysis 2, 1--21 (1994).

\bibitem{carr} T. Carroll, X. Massaneda, J. Ortega-Cerda, An enhanced uncertainty principle for the Vaserstein distance, Bull. London Math. Soc. 52, 1158--1173 (2020).

\bibitem{chung} F. Chung, Spectral Graph Theory, CBMS Regional Conference Series in Mathematics 92, American Mathematical Society (1997).

\bibitem{dahlberg} B. Dahlberg, Regularity Properties of Riesz Potentials, Indiana University Mathematics Journal 28, 257--268 (1979).

\bibitem{dick} J. Dick and F. Pillichshammer, Digital nets and sequences. Discrepancy theory and quasi-Monte Carlo integration, Cambridge University Press, Cambridge (2010).
  
\bibitem{ed} A. Edrei, Sur les d\'eterminants r\'ecurrents et les singularit\'es d’une fonction donn\'ee par son d\'eveloppement de Taylor, Compositio Mathematica 7, 20--88 (1940).

\bibitem{erd1}  P. Erd\H{o}s and P. Tur\'an, On a problem in the theory of uniform distribution. I. Nederl. Akad. Wetensch. 51, 1146--1154 (1948).

\bibitem{erd2}  P. Erd\H{o}s and P. Tur\'an, On a problem in the theory of uniform distribution. II. Nederl. Akad. Wetensch. 51, 1262--1269 (1948).

\bibitem{grig} A. Grigor'yan, Introduction to Analysis on Graphs, University Lecture Series 71, American Mathematical Society (2018).
  
\bibitem{hardin} D. Hardin and E. Saff, Discretizing manifolds via minimum energy points, Not. Amer. Math. Soc. (2004).
  
\bibitem{hu} P. Hu and W. C. Lau, A Survey and Taxonomy of Graph Sampling, arXiv:1308.5865 (2013).

\bibitem{kant} L.V. Kantorovich, On the Translocation of Masses. J Math Sci 133, 1381--1382 (2006).

\bibitem{koksma} J. Koksma, Een algemeene stelling uit de theorie der gelijkmatige verdeeling modulo 1, Mathematica B (Zutphen) 11, 7--11 (1942/43).

\bibitem{leja} F. Leja, Sur certaines suites li\'ees aux ensembles plans et leur application `a la repr\'esentation conforme, Annales Polonici Mathematici 4, 8--13 (1957).

\bibitem{linderman} G. Linderman and S. Steinerberger, Numerical Integration on Graphs: where to sample and how to weigh, Math. Comp. 89 (2020).

\bibitem{lopez} G. A. L\'opez and E. B. Saff. Asymptotics of greedy energy points. Math. Comp. 79 (272), 2287--2316 (2010).

\bibitem{num15} S. De Marchi and G. Elefante, Quasi-Monte Carlo integration on manifolds with mapped low-discrepancy points and greedy minimal Riesz s-energy points, Applied Numerical Mathematics 127, 110--124 (2018).

\bibitem{marzo} J. Marzo and A. Mas, Discrepancy of Minimal Riesz Energy Points, arXiv:1907.04814 (2019).

\bibitem{pausinger} F. Pausinger, Greedy energy minimization can count in binary: point charges and the van der Corput sequence, Annali di Matematica 200, 165--186 (2021).

\bibitem{pes1} I. Pesenson, Sampling in Paley-Wiener spaces on combinatorial graphs, Trans. Amer. Math. Soc. 360, 5603--5627 (2010).

\bibitem{pes2} I. Pesenson and M. Pesenson, Sampling, filtering and sparse approximations on combinatorial graphs, J. Fourier Anal. Appl. 16, 921--942 (2010).

\bibitem{peyre} G. Peyr\'{e} and M. Cuturi, Computational optimal transport, Foundations and Trends in Machine Learning 11, no. 5-6, 355--607 (2019).

\bibitem{pey}  R. Peyre, Comparison between $W_2$ distance and $\dot H^{-1}$ norm, and Localization of Wasserstein distance, ESAIM: COCV 24, 1489--1501 (2018).

\bibitem{num2} L. Reichel, Newton Interpolation at Leja points, BIT 30, 332--346 (1990).

\bibitem{num3} L. Reichel, The Application of Leja Points to Richardson Iteration and Polynomial Preconditioning, Linear Algebra and its Applications 154--156, 389--414 (1991).

\bibitem{saff} E. Saff and V. Totik, Logarithmic potentials with external fields, Springer (2013).

\bibitem{schmidt} W. M. Schmidt, Irregularities of distribution. VII. Acta Arith. 21, 45--50 (1972).

\bibitem{shu} D. Shuman, S. Narang, P. Frossard, A. Ortega and P. Vandergheynst, The emerging field of signal processing on graphs: Extending high-dimensional data analysis to networks and other irregular domains, IEEE Signal Process. Mag. 30(3), 83--98 (2013).

\bibitem{stein} S. Steinerberger, Wasserstein Distance, Fourier Series and Applications, arXiv:1803.08011 (2018).

\bibitem{steini} S. Steinerberger, Dynamically Defined Sequences with Small Discrepancy,  arXiv:1902.03269 (2019).

\bibitem{steinii} S. Steinerberger, Polynomials with Zeros on the Unit Circle: Regularity of Leja Sequences, arXiv:2006.10708 (2020).

\bibitem{thomson} J. J. Thomson, On the Structure of the Atom: an Investigation of the Stability and Periods of Oscillation of a number of Corpuscles arranged at equal intervals around the Circumference of a Circle; with Application of
the Results to the Theory of Atomic Structure, Philosophical Magazine Series 6, Volume 7, Number 39, 237--265 (1904).

\bibitem{wasser} L. N. Vasershtein, Markov processes on a countable product space describing large systems of
automata, Problemy Peredavci Informacii, 3, 64--72 (1969).

\bibitem{villani} C. Villani, Topics in Optimal Transportation, Graduate Studies in Mathematics, American Mathematical Society (2003).

\end{thebibliography}
\end{document}